\documentclass[12pt]{article}
\usepackage{amscd, amssymb, latexsym, amsmath}

\newtheorem{thm}{Theorem}

\newtheorem{que}{Question}

\newenvironment{proof}
{\noindent {\em Proof}} {\hfill $\Box$}

\numberwithin{thm}{section} \numberwithin{cor}{section}
\numberwithin{pro}{section} \numberwithin{lem}{section}
\numberwithin{dfn}{section}
\numberwithin{rem}{section} \numberwithin{equation}{section}

\newcommand{\R}{\mathbb R}

\begin{document}
\title{Mean Curvature Flows in Higher Codimension}
\author{Mu-Tao Wang }

\date{March 11, 2002}
\maketitle \centerline{email:\, mtwang@math.columbia.edu
\footnote{The author is supported in part by N.S.F. grant DMS
0104163.}}
\begin{abstract}
  The mean curvature flow is an evolution
process under which a submanifold deforms in the direction of its
mean curvature vector.  The hypersurface case has been much
studied since the eighties. Recently, several theorems on
regularity, global existence and convergence of the flow in
various ambient spaces and codimensions were proved. We shall
explain the results obtained as well as the techniques involved.
The potential applications in symplectic topology and mirror
symmetry will also be discussed.
\end{abstract}

\section{Introduction}

Let $M$ be a Riemannian manifold and $F:\Sigma \mapsto M$ an
isometric immersion of a smooth compact submanifold.  The second
fundamental form $A$ of $\Sigma$ is defined by

\[A:T\Sigma\times T\Sigma \mapsto N\Sigma, \]
\[A(X,Y)=(\nabla_X^M Y)^{\perp}\,\,\text{where}\, X, Y \in T\Sigma\]

Here $\nabla^M$ is the Levi-Civita connection on $M$ and
$(\cdot)^\perp$ denotes the projection onto the normal bundle
$N\Sigma$.

The mean curvature vector $H$  is defined by
\[H=Tr_g A\in N\Sigma\]
where the trace $Tr_g$ is taken with respect to the induced metric
$g$ on $\Sigma$.

When $M=\R^N$, if $x^1,\cdots ,x^n$ denote a local coordinate
system on $\Sigma$ then the second fundamental form can be
represented by $(\frac{\partial^2 F}{\partial x^i
\partial x^j})^\perp$ and the mean curvature vector is $H=(g^{ij}
\frac{\partial^2 F}{\partial x^i \partial x^j})^\perp$, where
$g^{ij}$ is the inverse matrix to $g_{ij}=< \frac{\partial
F}{\partial x^i}, \frac{\partial F}{\partial x^j}>$, the first
fundamental form. Here $\frac{\partial^2 F}{\partial x^i \partial
x^j}$ is considered as a vector in $\R^N$. In this case, it is not
hard to check

\begin{equation}\label{heat}
\,H=\Delta_\Sigma F
\end{equation}
 where $\Delta_\Sigma$ is the Laplace operator of the induced
metric on $\Sigma$.

The importance of the mean curvature vector lies in the first
variation formula of area.

\[\delta_V(area (\Sigma))=-\int_\Sigma H\cdot V\] for any variation
field $V$. Thus $H$ is the normal vector field on $\Sigma$ that
points to the direction in which the area decreases most rapidly.
$\Sigma$ is called a minimal submanifold if $H$ vanishes
identically.

The mean curvature flow of $F:\Sigma \mapsto M$ is a family of
immersions $F:\Sigma \times [0, \epsilon) \mapsto M$ parametrized
by $t$  that satisfies
\begin{equation}\label{eq}
\begin{split}
&\frac{d}{dt} F_t(x)=H(x,t)\\
&F_0=F
\end{split}
\end{equation}
where $H(x,t)$ is the mean curvature vector of $F_t(\Sigma)$ at
$F_t(x)$.

 This should be considered as the heat equation for submanifolds in view of equation (\ref{heat}). A submanifold
tends to find its optimal shape inside the ambient manifold.

 If we assume $M=\R^N$, in terms of coordinate $x^1,\cdots ,x^n$
on $\Sigma$, the mean curvature flow is the solution to the
following system of parabolic equations
\[F=F^A(x^1,\cdots, x^n, t),\,\, A=1,\cdots N\]
\[\frac{\partial F^A}{\partial t}
=\sum_{i, j, B}g^{ij}\,P^A_B\,\frac{\partial^2 F^B}{\partial
x^i\partial x^j}, \,\, A=1,\cdots N \] where $P^A_B
=\delta^A_B-g^{kl}\frac{\partial F^A} {\partial x^k}
\frac{\partial F^B}{\partial x^l}$ is the projection operator to
the normal direction.

 The equation (\ref{eq}) is a quasi-linear parabolic system and short time existence is
guaranteed  when the initial submanifold $\Sigma$ is compact and
smooth \cite{lsu}.

In general, the mean curvature flow fails to exist after a finite
time. The singularity is completely characterized by the blow up
of the second fundamental form. Namely, singularity at $t_0$ if
and only if $\sup_{\Sigma_t} |A|^2\rightarrow \infty$ as
$t\rightarrow t_0$. See for example \cite{hu1} for the
hypersurface case.

 The mean curvature flow has been studied by
 various approaches. In this article, we shall concentrated on
 the approaches of classical partial differential equations and
 geometric measure theory. For the level set approach and
 numerical methods, please see \cite{as} and the reference
 therein.

There are many beautiful results in the hypersurface (codimension
one) case.

\begin{thm}(Huisken, 1984 ($N\geq 3$)\cite{hu1}, Gage-Hamilton, 1985
($N=2$)\cite{gh}) Any convex compact hypersurface in $\R^N$
contracts to a round point after finite time along the mean
curvature flow.
\end{thm}

\begin{thm}(Grayson, 1987 \cite{gr})
Any embedded closed curve in $\R^2$ contracts to a round point
after finite time along the curvature flow.
\end{thm}

 In codimension one case,
$H$ is essentially a scalar function and $H>0$ is preserved along
the flow \cite{hu1}. As a contrast, in higher codimension $H$ is a
genuine vector and we do not know how to control the direction of
$H$. There are relatively very few results in the higher
codimension case, see  \cite{al}, \cite{as} and \cite{sm}.

 We shall discuss some new results about mean curvature flows in higher codimension in this article.  The
guideline is to identify positive quantities preserved along the
flow.

The author would like to thank Gerhard Huisken, Richard Hamilton,
Richard Schoen, Leon Simon, Brian White and Shing-Tung Yau for
valuable suggestions and discussions .

\section{Applications in calibrated geometry}

Let $M$ be an $N$-dimensional Riemannian manifold and $\alpha$ an
$n$-form on $M$. For any $p\in M$ and $S$ any $n$-dimensional
oriented subspace of $ T_p M$. Represent $S=e_1 \wedge \cdots
\wedge e_n$, where $e_1, \cdots e_n$ is an oriented orthonormal
basis for $S$. Define the \textit{comass} of $\alpha$ , $|\alpha|$
by

\[ |\alpha|(p)=\sup_{S\subset T_pM}
\alpha(e_1 \wedge \cdots \wedge e_n)
\]
Following Harvey and Lawson \cite{hl},
  a closed form $\alpha$ with $|\alpha|(p)=1$ at each $p\in M$
  is called a \textit{calibrating} form.

An oriented closed $n$-dimensional submanifold $\Sigma$ of $M$ is
\textit{calibrated} by $\alpha$ if $\alpha(T_p \Sigma)=1$ for any
$p\in \Sigma$, or $\alpha|_\Sigma=vol|_\Sigma$. A fundamental fact
in calibrated geometry is: a calibrated submanifold minimizes area
in its homology class. This follows from Stoke's Theorem:

Let $\Sigma'$ be any submanifold with $[\Sigma']=[\Sigma]$, then
\[
\int_{\Sigma'} vol|_{\Sigma'} \geq  \int_{\Sigma'}
\alpha|_{\Sigma'} =\int_\Sigma \alpha|_\Sigma =\int_\Sigma
vol|_\Sigma
 \]

We are interested in the following two classes of calibrated
submanifolds.

(1) Let $(M, \omega)$ be a K\"ahler manifold and
$\alpha=\frac{1}{n!}\omega^n$. Any $2n$ dimensional complex
submanifold is calibrated by $\alpha$.

(2) Let $(M,\omega, \Omega)$ be Calabi-Yau of complex dimension
$m$ and $\Omega$ is the parallel holomorphic $(m, 0)$ form.
$\alpha= Re \Omega$ is then a calibrating form. A Lagrangian
submanifold calibrated by $\alpha$ is called a \textit{special
Lagrangian} submanifold.

Define the function $*\alpha(p)=\alpha (T_p\Sigma)$ on $\Sigma$,
we may use $*\alpha$ to measure how far $\Sigma$ is away from
being calibrated. On a calibrated submanifold $*\alpha \equiv 1$.
It turns out the condition $*\alpha >0$ can be used to rule out a
certain type of singularity.

 There are two types of finite time singularity depending on the blow-up rate of $|A|$. Denote by $t_0$  the
blow up time, then $\sup_{\Sigma_t} |A|^2\rightarrow \infty $ as
$t\rightarrow t_0$. The singularity is said to be
\textit{fast-forming} (type I) if there exists a $C>0$ such that

\[\sup_t |A|^2\leq \frac{C}{t_0-t}\]

Otherwise, the singularity is called \textit{slow-forming} (type
II). For embedded curve on the plane, only type I singularity
occurs.

\begin{thm}\cite{mu1}
Let $(M^4, \omega)$ be a K\"ahler-Einstein four-manifold, then a
symplectic surface, i.e. $*\omega>0$ remains symplectic along the
mean curvature flow and the flow does not develop any type I
singularity.
\end{thm}

The results in Theorem 2.1 were obtained in the summer of 1999 and
announced in February 2000 at Stanford's differential geometry
seminar. Theorem 2.1 was also proved by Chen-Tian \cite{ct} and
Chen-Li \cite{cl}.

When $M$ is a Calabi-Yau manifold of arbitrary dimension, we prove
the following theorem.

\begin{thm}\cite{mu1} Let $(M, \omega, \Omega)$ be a Calabi-Yau manifold, then
a Lagrangian submanifold with $*Re \Omega>0$ remains Lagrangian
and $*Re\, \Omega>0 $ along the mean curvature flow and the flow
does not develop any type I singularity.
\end{thm}

That being Lagrangian is preserved along the mean curvature flow
in K\"ahler-Einstein manifolds was  proved by Smoczyk in
\cite{sm}.

\vskip 10pt
\begin{proof}. In \cite{mu1}( see page 324, remark 5.1), we sketch a proof of this theorem.
For a Lagrangian submanifold of Calabi-Yau manifolds, the
fundamental equations are, (see for example \cite{sm2} or
\cite{ty})

\begin{equation}\label{phase}
\Omega|_\Sigma=e^{i\theta} vol|_\Sigma
\end{equation}

\[H=J(\nabla \theta)\]


Once we know being Lagrangian is preserved, then $\theta$
satisfies the heat equation

\[\frac{d}{dt}\theta=\Delta \theta\]

By equation (\ref{phase}), $*Re\Omega=\cos\theta$. A
straightforward calculation using $|H|^2=|\nabla \theta|^2$ shows

\[\frac{d}{dt}*Re\Omega=\Delta *Re\Omega+|H|^2*Re\Omega\]

The same argument in Proposition 5.2 of \cite{mu1} shows on a type
I blow up limit $\Sigma_\infty$, we will have $|H|^2=0$. Since any
type I blow up limit is smooth and satisfies $F^\perp=H$
\cite{hu2}, $\Sigma_\infty$ must be a flat space. White's
regularity theorem \cite{w} shows there is no type I singularity.

\end{proof}

\section{Applications in mapping deformations }

Let $f:\Sigma_1\mapsto \Sigma_2$ be a smooth map between compact
Riemannian manifolds. The volume form $\omega_1$ of $\Sigma_1$
extends to a parallel calibrating form on the product space
$\Sigma_1\times \Sigma_2$. Let $\Sigma$ be the graph of $f$ as a
submanifold of $\Sigma_1\times \Sigma_2$. On $\Sigma$, $*\omega_1
=Jac(\pi_1|_\Sigma)$ is  the Jacobian of the projection
$\pi_1:\Sigma_1\times \Sigma_2\mapsto \Sigma_1$ restricting to
$\Sigma$. Any  submanifold $\Sigma'$ of $ \Sigma_1\times\Sigma_2$
is a locally a graph over $\Sigma_1$ if $*\omega_1>0$ on $\Sigma'$
by the inverse function theorem.

We shall evolve $\Sigma $ by the mean curvature flow in
$\Sigma_1\times \Sigma_2$ .  If $*\omega_1>0$ is preserved along
the flow then each $\Sigma_t$ is a graph over $\Sigma_1$ and thus
the flow gives a deformation $f_t$ of the original map $f$. It
turns out a stronger inequality is preserved.

\begin{thm} \cite{mu1}
If a smooth map $f:S^2\mapsto S^2$ satisfies $*\omega_1 >
|*\omega_2|$ on the graph of $f$, then the inequality remains true
along the mean curvature flow, the flow exists smoothly for all
time and
 $f_t$ converges to a
constant map.
\end{thm}

 The assumption is the same
as $Jac(\pi_1|_\Sigma)>|Jac(\pi_2|_\Sigma)|$. In other words, if
we see more of $\Sigma$ from $\Sigma_1$ than from $\Sigma_2$, then
$\Sigma$ converges to some $\Sigma_1\times \{p\}$ eventually. This
is a natural geometric assumption and we believe such assumption
is necessary for higher codimension mean curvature flow.

This theorem is generalized to arbitrary dimension and codimension
in \cite{mu3} under a slightly stronger assumption.

\begin{thm}\cite{mu3}
Let $f:S^n\mapsto S^m$ be a smooth map. If
$*\omega_1>\frac{1}{\sqrt{2}} $ on the graph of $f$, then the mean
curvature flow of the graph of $f$ in $S^n\times S^m$ exists for
all time, remains a graph, and converges smoothly to the graph of
a constant map  at infinity.
\end{thm}

This theorem is true under various curvature assumptions, please
see \cite{mu3} for the more general version. $*\omega_1$ should be
considered as the inner product of the tangent space of $\Sigma$
and the tangent space of $S_n$. The condition
$*\omega_1>\frac{1}{\sqrt{2}}$ guarantees $T\Sigma$ is closer to
$TS^n$ than to any other competing directions.

When $*\omega_1=*\omega_2$, we proved the following theorem.
\begin{thm} \cite{mu2}
Let $f:S^2\mapsto S^2$ be a smooth map such that
$*\omega_1=*\omega_2>0$ on the graph of $f$, then the equality is
preserved along the mean curvature flow and $f_t$ converges to an
isometry of $S^2$.
\end{thm}
The condition translates to $f^*\omega_2=\omega_1$, or $f$ is an
area-preserving diffeomorphism (or symplectomorphism). Recall the
harmonic heat flow of Eells-Sampson considers the deformation
 of a map $f:M\mapsto N$ along the gradient flow of the energy
 functional. When the sectional curvature of the target $N$ is non-positive, the flow exists for all time and
 $f_t$ converges to a harmonic map as $t\mapsto \infty$. For
 maps of nonzero degree between two-spheres $f:S^2\mapsto S^2$, singularities do occur in the harmonic heat flow
 even after finite time. It is quite surprising that the mean
 curvature deformation exists for all time and converges.

 $f_t$ indeed
provides a path in the  diffeomorphism (symplectomorphism) group
of $S^2$.

\begin{thm}\cite{mu2}
Any area preserving diffeomorphism of two-sphere deforms to an
isometry through area preserving diffeomorphisms along the mean
curvature flow.
\end{thm}

 For a Riemann surface $\Sigma $ with positive genus, the same
result holds \cite{mu2} when $\Sigma$ has hyperbolic metric and
the map $f:\Sigma \mapsto \Sigma$ is homotopic to identity.

\section{Proof of Theorem 3.1}

We shall explain the techniques involved in the  proof of Theorem
3.1 as it is the first complete solution to a higher codimension
mean curvature flow.
\subsection{Maximum principle}

The maximum principle of parabolic systems  developed by R.
Hamilton \cite{ha} plays an important role in the study of
geometric evolution equations. The first step is to use maximum
principle to show the inequalities $*\omega_1+*\omega_2>0$ and
$*\omega_1-*\omega_2>0$ are preserved along the flow.

In fact,  if we denote the singular values of ${df}$ by
$\lambda_1$ and $\lambda_2$, then
\[*\omega_1=\frac{1}{\sqrt{(1+\lambda_1^2)(1+\lambda_2^2)}}\]
and
\[*\omega_2=\frac{\lambda_1
\lambda_2}{\sqrt{(1+\lambda_1^2)(1+\lambda_2^2)}}\]

Let $\eta_1=*\omega_1+*\omega_2$ and $\eta_2=*\omega_1-*\omega_2$,
then $0<\eta_1, \eta_2\leq 1$.  The following equations are
derived in \cite{mu1} (see \cite{mu3} for general parallel forms).

\begin{equation}\label{eta_1}
 \frac{d}{dt}\eta_1=\Delta
\eta_1+\eta_1|A_1|^2+\eta_1(1-\eta_1^2)
\end{equation}

\begin{equation}\label{eta_2}
 \frac{d}{dt}\eta_2=\Delta
\eta_2+\eta_2|A_2|^2+\eta_2(1-\eta_2^2)
\end{equation}
 where $A_1$
and $A_2$ are part of the second fundamental form with
$|A_1|^2+|A_2|^2=2|A|^2$.

 The assumption of  Theorem 3.1 implies $\eta_1, \eta_2>0$ initially.
 By maximum principle of parabolic equations, $\min_{\Sigma_t} \eta_i$ is nondecreasing.
 This guarantees $*\omega_1>|*\omega_2|$ is preserved.  Adding
equations (\ref{eta_1}) and (\ref{eta_1}), we get for
$\mu=*\omega_1$,

\begin{equation}\label{mu}
 \frac{d}{dt}\mu\geq \Delta\mu+c|A|^2\end{equation}
 where $c>0$ is $\min\{\eta_1, \eta_2\}$ at $t=0$.

\subsection{Blow-up analysis}

The blow-up analysis is used in proving long time existence of the
flow.
 First let us recall the blow-up analysis for minimal surfaces. To study a possible
 singularity $x_0$, we blow up  the minimal surface $\Sigma^n$ at $x_0$ by
  $B_\lambda:x \mapsto \lambda(x-x_0)$, $\lambda>0$. Any limit as $\lambda\rightarrow \infty$ is
  still
minimal since the minimal surface equation is invariant under the
scaling. It must be a cone as a consequence of the monotonicity
formula. A minimal cone is rigid in the following sense:  if it is
close enough to a plane, then it must be a plane. The closeness is
measured by the density function.

\[\Theta(x_0)=\lim_{r\rightarrow 0} \Theta(x_0, r)=\lim_{r\rightarrow 0} \frac{area(B(x_0, r)\cap
\Sigma)}{\omega^n r^n}\] where $\omega^n$ is the area of an
$n$-dimensional unit ball.  The monotonicity formula in minimal
surface theory says $\Theta(x_0, r)$ is non-increasing  as $r $
approaches $0$, in particular the limit exists. Allard's
regularity theorem \cite{all} then asserts there exists an
$\epsilon>0$ such that if $\Theta(x_0) <1+\epsilon$, then $x_0$ is
a regular point. We refer to Simon's book \cite{gmt} for minimal
surface theory.

 For the mean curvature flow, we consider the total space time as a submanifold
in $M\times \R$ and use parabolic blow up at a space time point
$(x_0, t_0)$. The limit is still a mean curvature flow  and  the
monotonicity formula of Huisken implies a time slice satisfies
$F^\perp=H$, or the limit flow is self-similar. For a smooth
point, we obtain the stationary flow of a plane.

 The density function is now replaced by the integral of the backward
heat kernel. To be more precise, we isometrically embed $M$ into
$\R^N$.  For any $\lambda
>1$, the parabolic blow up $D_\lambda$ at $(x_0, t_0)$ is defined
by

\begin{equation}
\begin{split}
D_\lambda :\,\R^N \times[0, t_0) &\rightarrow
\R^N \times [-\lambda^2 t_0, 0)\\
(x,t)&\rightarrow (\lambda(x-x_0), \lambda^2 (t-t_0))
\end{split}
\end{equation}

 The (n-dimensional) backward heat kernel $\rho_{x_0, t_0}$ at
$(x_0, t_0)$ is

\begin{equation}
\rho_{x_0, t_0}(x,t)=\frac{1}{(4\pi(t_0-t))^{n\over 2}} \exp
(\frac{-|x-x_0|^2}{4(t_0-t)})
\end{equation}

 Notice that the integral $\int \rho_{x_0, t_0} d\mu_t$ is invariant
under the parabolic blow up, where $d\mu_t=\sqrt{\det
g_{ij}(F_t)}d\mu$ is the pull back volume form by $F_t$.

 The monotonicity formula of Huisken \cite{hu2} says for $t<t_0$

\[\frac{d}{dt}\int \rho_{x_0, t_0} d\mu_t
\leq 0\] so the limit as $t\rightarrow t_0$ exists. This formula
holds only for mean curvature flows in Euclidean spaces. For a
general ambient manifold, a modification to take care of curvature
terms is necessary, see \cite{w2} or \cite{il1}. The analogue of
Allard's regularity theorem in mean curvature flow is the
following theorem of White.

\begin{thm}\cite{w}
There is an $\epsilon >0$ such that if
\[\lim_{t\rightarrow t_0}\int \,\rho_{x_0, t_0}
d\mu_t <1+\epsilon\] then $(x_0, t_0)$ is a regular point.
\end{thm}

In the proof of Theorem 3.1,  the equation (\ref{mu}) helps us
find a subsequence $t_i\rightarrow t_0$ and
 blow up rate $\lambda_i\rightarrow \infty$ such that the $L^2$ norm of the second
  fundamental form $\int |A|^2$ of  $\lambda_i
\Sigma_{t_i}$ approaches zero. The limit is thus a plane  and
$\int \rho_{x_0, t_0} d\mu_{t_i}\rightarrow 1$ as $t_i \rightarrow
\infty$. Monotonicity formula implies $\lim_{t\rightarrow t_0}\int
\,\rho_{x_0, t_0} d\mu_t =1$. By White's theorem again it can be
concluded that $(x_0, t_0)$ is a regular point.

\subsection{Curvature estimate}

 From the calculation in
\cite{mu1}, the norm of the second fundamental form $|A|^2$
satisfies

\[
\begin{split}
\frac{d}{dt} |A|^2
&\leq \Delta |A|^2 -2|\nabla A|^2+K_1|A|^4+K_2 |A|^2\\
\end{split}
\]
where $K_1, K_2$ are constants depending on the curvature and the
covariant derivatives of curvature  of the ambient space .

 In general, $|A|^2 $ blows up in finite time because of the
$|A|^4$ term. However, the equation of $\mu$ (\ref{mu}) helps to
control $|A|^2$ in the proof of Theorem 3.1.

 By equation (\ref{mu}) for any $k$, $0<k<1$,
 we have
\[\frac{d}{dt}(\mu-k)\geq\Delta (\mu-k)+c\frac{\mu}{\mu-k} (\mu-k)|A|^2
\]

where we use $0<\mu \leq 1$.

 If $\min_{\Sigma_t} \mu$ is very close to one when $t$ is large, we may choose $k$ close to $1$ so that
$\mu-k>0$ is preserved after some $t_1$ and $\frac{\mu}{\mu-k}$ is
large. The quantity $g=\frac{|A|^2}{\mu-k}$ after time $t_1$ then
satisfies
\[
\begin{split}
\frac{d}{dt}g \leq \Delta g +V\cdot \nabla g-K_3 g^2+K_4 g
\end{split}
\]
with $K_3>0$. The maximum principle shows $g$ is uniformly bounded
for $t>t_1$ .

By equations (\ref{eta_1}) and (\ref{eta_2}) and a comparison
argument, we see $\min_{\Sigma_t} \eta_1 \rightarrow 1$ and
$\min_{\Sigma_t} \eta_2 \rightarrow 1$ as $t\rightarrow \infty$.
In particular , $\min_{\Sigma_t} \mu \rightarrow 1$. The
assumption on $\mu$ is true when $t$ is large enough and thus
$|A|^2$ is uniformly bounded. Integrate equation (\ref{mu}) over
space and time shows $\int_{\Sigma_t}|A|^2\rightarrow 0$ and the
sub mean value inequality in \cite{il1} shows $\sup
|A|^2\rightarrow 0$. The last step is to apply Simon's \cite{si}
general convergence theorem for gradient flows.

\section{Related problems}

Recall \cite{sm} that being Lagrangian is preserved along the mean
curvature flow in K\"ahler-Einstein manifolds. The graph of a
symplectomorphism of a K\"ahler-Einstein manifold $M$ is a
Lagrangian submanifold in the product space $(M\times M,
\omega_1-\omega_2)$. The following question is thus a natural
generalization of Theorem 3.3.

\begin{que} Can one prove the long time existence and convergence
of  mean curvature flows of  symplectomorphisms of
K\"ahler-Einstein manifolds?
\end{que}

 Theorem 3.3 and the corresponding theorems for higher genus Riemann surfaces implies any Lagrangian
graph is Lagrangian isotopic to a minimal Lagrangian graph along
the mean curvature flow. This is related to the following
conjecture due to Thomas and Yau \cite{ty}.

\begin{que}
Can one prove the long time existence and convergence of mean
curvature flow of a stable graded Lagrangian submanifold in a
Calabi-Yau manifold?
\end{que}

A notion of stability for Lagrangian submanifolds was formulated
in \cite{ty} in terms of the range of the phase function $\theta$.

\[\Omega|_\Sigma=e^{i\theta}vol|_\Sigma\]
on any Lagrangian submanifold of Calabi-Yau manifold. Theorem 2.2
implies if $-\frac{\pi}{2}<\theta<\frac{\pi}{2}$ on a Lagrangian
submanifold, then the mean curvature flow does not develop any
type I singularity. How to exclude or perturb away type II
singularities seems a very interesting yet hard problem.


\begin{thebibliography}{99}

\bibitem{all} W. K. Allard, \textit{On the first variation of a varifold.} Ann.
of Math. (2) 95 (1972), 417--491.
\bibitem{al} S. Altschuler, \textit{ Singularities of the curve shrinking flow for
space curves.} J. Differential Geom. 34 (1991), no. 2, 491--514.

\bibitem{as} L. Ambrosio and H. M. Soner, \textit{  Level set approach to mean
curvature flow in arbitrary codimension.} J. Differential Geom. 43
(1996), no. 4, 693--737.

\bibitem{cl} J. Chen and J. Li, \textit{ Mean curvature flow of surface in
$4$-manifolds.} Adv. Math. 163 (2001), no. 2, 287--309.

\bibitem{ct} J. Chen and G. Tian, \textit{ Moving symplectic curves in
Kähler-Einstein surfaces.} Acta Math. Sin. (Engl. Ser.) 16 (2000),
no. 4, 541--548.

\bibitem{eh} K. Ecker and G. Huisken \textit{Mean curvature evolution of entire
graphs}, Ann. of Math. (2) 130 (1989), no. 3, 453--471.

\bibitem{eh2} K. Ecker and G. Huisken,
\textit{ Interior estimates for hypersurfaces moving by mean
curvature.}, Invent. Math. 105 (1991), no. 3, 547--569.

\bibitem{es} J. Eells and J. H. Sampson, \textit{Harmonic mappings of
Riemannian manifolds.} Amer. J. Math. 86 1964 109--160.

\bibitem{es2}C. J. Earle and J. Eells, \textit{ The diffeomorphism group of a compact
Riemann surface.} Bull. Amer. Math. Soc. 73 1967 557--559.

\bibitem{gh} M. Gage and R. Hamilton, \textit{The heat equation shrinking convex
plane curves.}, J. Differential Geom. 23 (1986), no. 1, 69--96.

\bibitem{gr} M. Grayson,\textit{ The heat equation shrinks embedded plane
curves to round points.} J. Differential Geom. 26 (1987), no. 2,
285--314.
\bibitem{ha} R. Hamilton, \textit{ Four-manifolds with positive curvature operator}. J.
Differential Geom. 24 (1986), no. 2, 153--179.



\bibitem{hl} R. Harvey and H. B. Lawson, \textit{ Calibrated
geometries.} Acta Math. 148 (1982), 47--157.

\bibitem{hu1} G. Huisken, \textit{Flow by mean curvature of convex surfaces into
spheres.} J. Differential Geom. 20 (1984), no. 1, 237--266.

\bibitem{hu2} G. Huisken, \textit{Asymptotic
behavior for singularities of the mean curvature flow}, J.
Differential Geom. \textbf{31} (1990), no. 1, 285--299.




\bibitem{il1} T. Ilmanen, \textit{Singularities
of mean curvature flow of surfaces }, preprint , 1997.

\bibitem{il2} T. Ilmanen, \textit
{Elliptic regularization and partial regularity for motion by mean
curvature},
 Mem. Amer. Math. Soc. 108 (1994), no. 520,

\bibitem{lsu} O. A. Ladyzenskaja, V. A. Solonnikov and N. N.
Uralceva \textit{Linear and quasilinear equations of parabolic
type}, Transl. Amer. Math. Soc. 23 (1968).

\bibitem{si} L. Simon, \textit{Asymptotics for
 a class of nonlinear evolution equations,
 with applications to geometric problems.},
  Ann. of Math. (2) 118 (1983), no. 3, 525--571.

\bibitem{gmt} L. Simon, \textit{ Lectures on geometric measure theory.}, Proceedings
of the Centre for Mathematical Analysis, Australian National
University, 3. Australian National University, Centre for
Mathematical Analysis, Canberra, 1983. vii+272 pp.

\bibitem{sm}  K. Smoczyk, \textit{
A canonical way to deform a Lagrangian
 submanifold.}, preprint, dg-ga/9605005.


\bibitem{sm2} K. Smoczyk,\textit{ Harnack inequality for the Lagrangian mean
curvature flow.} Calc. Var. Partial Differential Equations 8
(1999), no. 3, 247--258.


\bibitem{ty}  R. P. Thomas, S.-T. Yau \textit{ Special Lagrangians, stable bundles and mean
curvature flow.}, math.DG/0104197.

\bibitem{mu1}  M-T. Wang, \textit{Mean Curvature
Flow of surfaces in Einstein Four-Manifolds }, J. Differential
Geom. \textbf{57} (2001), no. 2, 301-338.

\bibitem{mu2} M-T. Wang, \textit{Deforming area preserving
diffeomorphism of surfaces by mean curvature flow }, Math. Res.
Lett. 8 (2001), no.5-6, 651-662.


\bibitem{mu3}  M-T. Wang, \textit{Long-time existence and
convergence of graphic mean curvature flow in arbitrary
codimension }, to appear in Invent. Math.




\bibitem{w} B. White, \textit{A local
regularity theorem for classical mean curvature flow}, preprint,
2000.
\bibitem{w2}B. White, \textit{Stratification of minimal surfaces, mean curvature flows,
 and harmonic maps.} J. Reine Angew. Math. 488 (1997), 1--35.
\end{thebibliography}
\end{document}